\documentclass[12pt]{article}
\usepackage[colorlinks,
            linkcolor=blue,
            anchorcolor=green,
            citecolor=blue
            ]{hyperref}
\usepackage{mathrsfs,amssymb,amsfonts}
\usepackage{amsmath,amssymb,latexsym,color}
\usepackage[mathscr]{eucal}
\usepackage[numbers,sort&compress]{natbib}
\usepackage{graphicx}
\usepackage{lipsum}

\def\thefootnote{\fnsymbol{footnote}}
\newtheorem{thm}{Theorem}[section]

\newtheorem{lemma}[thm]{Lemma}
\newtheorem{cor}[thm]{Corollary}
\newtheorem{example}[thm]{Example}

\newtheorem{problem}[thm]{Problem}
\newtheorem{ques}[thm]{Question}
\newtheorem{remark}[thm]{Remark}
\newtheorem{Notation}[thm]{Notation}

\newcommand{\proof}{{\it Proof.\quad}}
\newcommand{\qed}{\hfill\Box\medskip}
\renewcommand{\thefootnote}{\arabic{footnote}}

\newcommand{\diam}{{\rm diam}}
\newcommand{\Cyc}{{\rm Cyc}}
\begin{document}

\renewcommand{\abovewithdelims}[2]
{\genfrac{[}{]}{0pt}{}{#1}{#2}}

\title{\bf Notes on the diameter of the complement of the power graph of a finite group}

\author{{Xuanlong Ma$^{a}$, Alireza Doostabadi$^{b}$, Kaishun Wang$^{c}$}
\\
{\small\em $^{a}$School of Science, Xi'an Shiyou University, Xi'an 710065, China}
\\
{\small\em $^{b}$Faculty of Science, University of Zabol,
Zabol, Iran}
\\
{\small\em $^{c}$Laboratory of Mathematics and Complex Systems (MOE),}
\\
{\small\em School of Mathematical Sciences, Beijing Normal University, Beijing 100875, China}
}

 \date{}
 \maketitle
\newcommand\blfootnote[1]{%
\begingroup
\renewcommand\thefootnote{}\footnote{#1}%
\addtocounter{footnote}{-1}%
\endgroup
}
\begin{abstract}
We determine the diameter of every connected component of the complement of the power graph and the enhanced power graph
of a finite group,
 which completely answers two questions by Peter J. Cameron.


\end{abstract}

\blfootnote{E-mail addresses: xuanlma@mail.bnu.edu.cn (X. Ma), aldoostabadi@uoz.ac.ir (A. Doostabadi), wangks@bnu.edu.cn (K. Wang)}

\section{Introduction and preliminaries}
All graphs considered in this paper are finite, undirected, with no loops and no multiple edges. Let $\Gamma$ be a graph.
The vertex set and the edge set of $\Gamma$ are denoted by $V(\Gamma)$ and $E(\Gamma)$, respectively.
The {\em complement graph} of $\Gamma$, denoted by $\overline{\Gamma}$, is the graph whose vertex set is $V(\Gamma)$ and whose edges are the
non-edges of $\Gamma$.
Let $x,y\in V(\Gamma)$.
The {\em  distance} between $x$ and $y$ in $\Gamma$  is the length
of a shortest path from $x$ to $y$ and is denoted by $d_{\Gamma}(x,y)$.
If the situation is unambiguous, we denote $d_{\Gamma}(x,y)$ simply by $d(x,y)$.
If two distinct vertices $x$
and $y$ of $\Gamma$ are adjacent, then we denote this by
$x \sim_{\Gamma} y$, or shortly by $x \sim y$.
Particularly, we denote by $x_1\sim x_2 \sim \cdots \sim x_n$
a path of $n$ vertices in $\Gamma$.
The {\em diameter} of $\Gamma$, denoted by $\diam(\Gamma)$, is the greatest distance between any two vertices.
Note that $\diam(\Gamma)=1$ if and only if $\Gamma$ is a complete graph.

Graphs associated with groups and other algebraic structures have been actively investigated, since they have valuable applications
(cf. \cite{KeR}) and are related to automata theory (cf. \cite{K}).
The {\em undirected power graph} $\mathcal{P}(G)$ of a finite group $G$ has vertex set $G$ and two distinct elements are adjacent if one is a power of the other.
The concepts of power graph
and undirected power graph were first introduced by Kelarev and Quinn \cite{KQ00} and Chakrabarty {\em et al.} \cite{CGS}, respectively.
In recent years, the study of power graphs has been growing, for example, Cameron, Manna and Mehatari \cite{Cam221} studied the finite groups whose power graphs are cographs; Zahirovi\'c \cite{Za21} explored the directed power graph of a torsion-free group determined by its power graph, and showed that any two torsion-free groups having isomorphic power graphs have isomorphic directed power graphs; Ma, Feng and Wang \cite{Ma21} investigated the Lambda number of the power graph of a finite group; Manna, Cameron and Mehatari \cite{MCam21} studied some forbidden subgraphs of power graphs of groups and gave a number of open problems. Also, see,
for example, \cite{CamJ,B2,Doo,Doo1}.
Let us refer to two surveys \cite{AKC,KSC21} for more information pertaining to the research results and open
problems on the power graphs of groups.

In order to measure how close the power graph is to the commuting graph, Aalipour  {\em et al.} \cite{acam} introduced the enhanced power graph of a group which lies in between.
Let $G$ be a finite group.
The {\em enhanced power graph} $\mathcal{P}_e(G)$ of $G$ is the graph whose vertex set is $G$, and two distinct vertices are adjacent if they generate a cyclic subgroup of $G$.
The enhanced power graph has also appeared in the literature under the name cyclic graph (cf. \cite{CM2}).
In recent years, the study of enhanced power graphs has received considerable attention. For example,
Bera and Bhuniya \cite{Bh18} showed that there is a one-to-one correspondence between the
maximal cliques in $\mathcal{P}_e(G)$ and the maximal cyclic subgroups of $G$. In 2020, Zahirovi\'c, Bo\v{s}njak and Madar\'asz \cite{Za20} showed that any isomorphism between undirected power graph of finite groups is an isomorphism between enhanced power graphs of these groups.
Ma and She \cite{MS} characterized the metric dimension of an enhanced power graph.
Bera, Dey and Mukherjee \cite{Bera} completely characterized the abelian groups such that their proper enhanced power graphs are connected, where the {\em proper enhanced power graph} of $G$ is the induced subgraph of $\mathcal{P}_e(G)$ obtained by deleting the identity element of $G$.
The reader is referred to the survey paper \cite{MKer} for a large number of results on enhanced power graphs of groups.

Every group considered in this paper is finite. We always use $e$ to denote the identity
element of the group under consideration.
Denote by $\mathbb{Z}_n$ the cyclic group of order $n$.
Let $G$ be a group.
A {\em maximal cyclic subgroup} of $G$ is a cyclic subgroup, which is not a proper subgroup of some cyclic
subgroup of $G$. The set of all maximal cyclic subgroups of $G$
is denoted by $\mathcal{M}(G)$. Note
that $|\mathcal{M}(G)|=1$ if and only if $G$ is cyclic.
For a subset $S$ of $G$, define
$$
\mathcal{M}_S:=\{M\in \mathcal{M}(G): S\subseteq M\}.
$$
If $S=\{s\}$, then we denote $\mathcal{M}_{\{s\}}$ simply by $\mathcal{M}_s$.

Recently, Cameron \cite{Cam22}
introduced aspects of various graphs whose vertex set is a group and whose edges
reflect group structure in some way, and proposed the following two questions.

\begin{ques}\label{ques1}{\rm (\cite[Question 19]{Cam22})}
What is the best possible upper bound for the diameter of non-trivial connected component of the complement of the power graph, and which groups attain the bound?
\end{ques}

\begin{ques}\label{ques2}{\rm (\cite[Question 20]{Cam22})}
Is it true that the complement of the enhanced power graph has just one connected component, apart from isolated vertices?
\end{ques}

Moveover, in \cite{Cam22}, Cameron proved the following result.

\begin{thm}\label{y-thm1}{\rm (\cite[Theorem 9.9]{Cam22})}
Let $G$ be a finite group which is not a cyclic $p$-group.
Then $\overline{\mathcal{P}(G)}$ has just one connected component, apart from isolated vertices.
\end{thm}

Note by \cite[Theorem 2.12]{CGS} that if $G$ is a cyclic $p$-group, then $\mathcal{P}(G)$ is complete, and so $\overline{\mathcal{P}(G)}$ has no edges, that is, every vertex is isolated. Thus, for $\overline{\mathcal{P}(G)}$, we always consider the finite groups $G$ which are not cyclic $p$-groups.
Moreover, in view of Theorem~\ref{y-thm1}, if $G$ is not a cyclic $p$-group, then all non-isolated vertices of $\overline{\mathcal{P}(G)}$ will induce a connected component, and we denote the connected component by $\overline{\mathcal{P}(G)^{\ast}}$.

Note that a finite group $G$ is cyclic if and only if
$\mathcal{P}_e(G)$ is complete (cf. \cite{acam,Bh18}). Thus, if $G$ is cyclic, then
$\overline{\mathcal{P}_e(G)}$ is an empty graph, that is, every vertex is isolated.
The {\em cycle} \cite{Ob92} of a group $G$, denoted by $\Cyc(G)$, is defined by
\begin{equation*}\label{peq1}
\Cyc(G)=\{g\in G: \langle g,x \rangle \text{ is cyclic for any $x\in G$}\}.
\end{equation*}
Clearly, a vertex $g$ in $\overline{\mathcal{P}_e(G)}$ is isolated if and only if $g\in \Cyc(G)$ (cf. \cite{MS}).
Denote by $\overline{\mathcal{P}_e(G)^{\ast}}$ the induced subgraph of $\overline{\mathcal{P}_e(G)}$ by all non-isolated vertices in $\overline{\mathcal{P}_e(G)}$. Namely,
$$
V(\overline{\mathcal{P}_e(G)^{\ast}})=G\setminus \Cyc(G).
$$

For $n\ge 3$,
the {\em generalized quaternion group} (also called {\em dicyclic group}) $Q_{2^n}$ of order $2^n$ has
a presentation
\begin{equation}\label{gq}
Q_{2^n}=\langle x,y: x^{2^{n-2}}=y^2, x^{2^{n-1}}=e, y^{-1}xy=x^{-1}\rangle.
\end{equation}
By \eqref{gq}, it is easy to check that
\begin{equation}\label{gq-1}
Q_{2^n}=\langle x\rangle\cup \{x^iy: 1\le i \le 2^{n-1}\}, ~~
o(x^iy)=4 \text{ for any $1\le i \le 2^{n-1}$}.
\end{equation}
We remark that $Q_{2^n}$ has the unique involution $y^2$ and
\begin{equation}\label{gq-2}
V(\overline{\mathcal{P}(Q_{2^n})^{\ast}})
=Q_{2^n}\setminus\{e,y^2\}
\end{equation}
by \cite[Proposition~4]{Cam10} or \cite[Theorem 9.1(a)]{Cam22}.

In order to state our main results, we first define a class of finite non-$p$-groups. A finite group $G$ is called a {\em $\Psi$-group} if the following hold:
\begin{itemize}
\item $|G|=p_1^{\alpha_1}p_2^{\alpha_2}\cdots p_t^{\alpha_t}$,
where $p_1,p_2,\ldots,p_t$ are distinct primes and $t\ge 2$;
  \item For any $1\le i \le t$, $G$ has a unique subgroup of order $p_i$;
  \item $G$ has an element of order
  $p_1^{\beta_1}p_2^{\alpha_2}\cdots p_t^{\alpha_t}$ such that $1\le \beta_{1}< \alpha_{1}$ and for any prime $p\ne p_1$, the Sylow $p$-subgroup of $G$ is unique.
\end{itemize}
By \cite[Theorem~5.4.10(ii)]{Gor} and the definition of a $\Psi$-group, we have the following remark.

\begin{remark}
Suppose that $G$ is a $\Psi$-group of order $p_1^{\alpha_1}p_2^{\alpha_2}\cdots p_t^{\alpha_t}$ and has an element of order $p_1^{\beta_1}p_2^{\alpha_2}\cdots p_t^{\alpha_t}$, where $t\ge 2$ and $1\le \beta_{1}< \alpha_{1}$.
  Then $G$ is isomorphic to one of the following:
$$
H \rtimes \mathbb{Z}_{p_1^{\alpha_1}},H\rtimes Q_{2^{n}},
H \times \mathbb{Z}_{p_1^{\alpha_1}},
H \times Q_{2^{n}},
$$
where
$H=\mathbb{Z}_{p_2^{\alpha_2}p_3^{\alpha_3}\cdots p_t^{\alpha_t}}$ and $n\ge 3$.
\end{remark}

We then define a class of non-cyclic groups.
A finite non-cyclic group $G$ is called a {\em $\Phi$-group} if there exist $x,y\in V(\overline{P_e(G)^{\ast}})$ such that the following hold:
\begin{itemize}
\item $\langle x,y\rangle$ is cyclic;
\item $\langle x\rangle\notin \mathcal{M}(G)$ and $\langle y\rangle\notin \mathcal{M}(G)$;
\item For any $M\in \mathcal{M}(G)\setminus \mathcal{M}_{\{x,y\}}$, either $x\in M$ or $y\in M$.
\end{itemize}

In this paper, we completely answer Questions~\ref{ques1} and \ref{ques2}.
Our main results are the following theorems.

\begin{thm}\label{mainthm}
Let $G$ be a finite group which is not a cyclic $p$-group.
Then
$$\diam(\overline{\mathcal{P}(G)^{\ast}})
=\left\{
                                  \begin{array}{ll}
                                  1, & \hbox{if $G\cong \mathbb{Z}_2^m$, where $m$ is a positive integer at least $2$;} \\
                                  3, & \hbox{if $G$ is a $\Psi$-group;} \\
                                  2, & \hbox{otherwise,}
                                  \end{array}
                                \right.
$$
\end{thm}

\begin{thm}\label{mainthm1}
Let $G$ be a finite non-cyclic group.
Then
$$\diam(\overline{\mathcal{P}_e(G)^{\ast}})
=\left\{
                                  \begin{array}{ll}
                                  1, & \hbox{if $G\cong \mathbb{Z}_2^m$, where $m$ is a positive integer at least $2$;} \\
                                  3, & \hbox{if $G$ is a $\Phi$-group;} \\
                                  2, & \hbox{otherwise,}
                                  \end{array}
                                \right.
$$
In particular, $\overline{\mathcal{P}_e(G)}$ has just one connected component, apart from isolated vertices.
\end{thm}

The next results are obtained
by applying Theorems~\ref{mainthm} and \ref{mainthm1} to
$p$-groups and cyclic groups.

\begin{cor}
Suppose that $G$ is a finite $p$-group which is non-cyclic. Then
$$\diam(\overline{\mathcal{P}_e(G)^{\ast}})
=\left\{
                                  \begin{array}{ll}
                                  1, & \hbox{if $G\cong \mathbb{Z}_2^m$, where $m$ is a positive integer at least $2$;} \\
                                  2, & \hbox{otherwise,}
                                  \end{array}
                                \right.
$$
\end{cor}

\begin{cor}
Suppose that $G$ is a finite cyclic group which is not a $p$-group. Then
$$\diam(\overline{\mathcal{P}(G)^{\ast}})
=\left\{
                                  \begin{array}{ll}
                                  2, & \hbox{if $|G|$ is a product of distinct primes;} \\
                                  3, & \hbox{otherwise,}
                                  \end{array}
                                \right.
$$
\end{cor}

\begin{cor}
Suppose that $G$ is a finite nilpotent group which is not a cyclic group. Then
$$\diam(\overline{\mathcal{P}(G)^{\ast}})
=\left\{
                                  \begin{array}{ll}
                                  1, & \hbox{if $G\cong \mathbb{Z}_2^m$, where $m$ is a positive integer at least $2$;} \\
                                  3, & \hbox{if $G\cong Q_{2^m}\times \mathbb{Z}_n$, where $m\ge 3$ and $n\ge 3$ with $2\nmid n$;} \\
                                  2, & \hbox{otherwise,}
                                  \end{array}
                                \right.
$$
\end{cor}

\section{Proof of Theorem~\ref{mainthm}}

In this section, we will prove Theorem~\ref{mainthm}.
We next give some lemmas before giving the proof of
Theorem~\ref{mainthm}.

Recall first the following elementary result.

\begin{lemma}\label{eresult}
{\rm (\cite[Theorem~5.4.10(ii)]{Gor})}\label{uniquep}
A $p$-group having a unique subgroup of order $p$ is either cyclic or generalized quaternion, where $p$ is a prime.
\end{lemma}

\begin{lemma}\label{p-subgroup}
Let $G$ be a finite group which is not a cyclic $p$-group, where $p$ is a prime.
For some divisor $p$ of $|G|$, if there exists two distinct subgroups of order $p$, then $\diam(\overline{\mathcal{P}(G)^{\ast}})\le 2$.
\end{lemma}
\proof
Let $\langle a\rangle$ and $\langle b\rangle$ be two distinct
subgroups of order $p$. Assume that $x$ and $y$ are non-adjacent in $\overline{\mathcal{P}(G)^{\ast}}$. It suffices to prove that $d(x,y)=2$. Note that in this case $\langle x\rangle\subseteq \langle y\rangle$ or $\langle y\rangle\subseteq \langle x\rangle$.
Without loss of generality, now let $\langle x\rangle\subseteq \langle y\rangle$. If $p\mid o(x)$, then it follows that one of $a$ and $b$ can not belong to $\langle x\rangle$, say $a\notin \langle x\rangle$, and so $x\sim a \sim y$ is a path, which implies $d(x,y)=2$, as desired. Suppose next that $p\nmid o(x)$. Note that one of $a$ and $b$ must not belong to $\langle y\rangle$, say $a\notin \langle y\rangle$. It follows that
$x\sim a \sim y$ is a path, which also implies $d(x,y)=2$, as desired.
$\qed$

\begin{lemma}\label{p-group}
Suppose that $G$ is a $p$-group which is non-cyclic. Then $$\diam(\overline{\mathcal{P}(G)^{\ast}})
=\left\{
                                  \begin{array}{ll}
                                  1, & \hbox{if $G\cong \mathbb{Z}_2^m$ for some positive integer $m$;}\\
                                  2, & \hbox{otherwise,}
                                  \end{array}
                                \right.
$$
\end{lemma}
\proof
Clearly, if every element of $G$ has order at most $2$, then $G\cong \mathbb{Z}_2^m$ for some positive integer $m$, and so
$\overline{\mathcal{P}(G)^{\ast}}$ is complete, which implies
$\diam(\overline{\mathcal{P}(G)^{\ast}})=1$. Thus, in the following, we may assume that $G$ has an element of order at least $3$. In view of Lemma~\ref{p-subgroup}, we may assume that $G$ has a unique subgroup of order $p$. It follows from Lemma~\ref{eresult} that $G$ is a generalized quaternion group.
Suppose that $x$ and $y$ are two non-adjacent vertices of $\overline{\mathcal{P}(G)^{\ast}}$.
Without loss of generality, let $\langle x\rangle\subseteq \langle y\rangle$.
By \eqref{gq-2}, it is easy to see that $4\mid o(x)$.
Combining now \eqref{gq} and \eqref{gq-1}, we have that there exists an element $z$ of order $4$ such that $z\notin \langle x\rangle$. It follows that $x\sim z\sim y$ is a path, and hence $d(x,y)=2$, as desired.
$\qed$

We next consider the finite groups that are not $p$-groups.
For a positive integer $n$, denote by $\pi(n)$ the set of all prime divisors of $n$.

\begin{lemma}\label{lem-0}
Let $|G|=p_1^{\alpha_1}p_2^{\alpha_2}\cdots p_t^{\alpha_t}$ where $t\ge 2$. Suppose that for any prime $p_i$, $G$ has a unique subgroup of order $p_i$ where $1\le i \le t$.
If $x$ and $y$ are distinct vertices of $\overline{P(G)^{\ast}}$ such that
$\langle x\rangle\subseteq \langle y\rangle$ and
$|G|/o(y)\ne p_i^{\beta_i}$ with $1\le \beta_i < \alpha_i$, then $d(x,y)=2$.
\end{lemma}
\proof
Note that $o(y)\ne |G|$, since $y\in V(\overline{P(G)^{\ast}})$.
If $\pi(o(y))\ne \pi(|G|)$,
taking $p\in \pi(|G|)\setminus \pi(o(y))$, we have that $G$ has an element $z$ of order $p$, and it follows that $x\sim z\sim y$ is a path, as desired.

Thus, in the following, we may assume that $\pi(o(y))=\pi(|G|)$.
Let
$$
o(y)=p_1^{\gamma_1}p_2^{\gamma_2}\cdots p_t^{\gamma_t},
~~1\le \gamma_i \le \alpha_i.
$$
Since $|G|/o(y)\ne p_i^{\beta_i}$ with $1\le \beta_i < \alpha_i$, there exist two distinct primes $p_i,p_j$ such that
$\gamma_i<\alpha_i$ and $\gamma_j<\alpha_j$. Note that one of $p_i$ and $p_j$ must not equal to $2$, without loss of generality, say $p_i\ne 2$.
Thus, Lemma~\ref{eresult} implies that $G$ has a cyclic Sylow $p_i$-subgroup of order $p_i^{\alpha_i}$, say $\langle w\rangle$. In this case, if $o(x)$ is not a power of $p_i$, then
it is clear that $x\sim w\sim y$ is a path, and so $d(x,y)=2$.
Now suppose that $o(x)=p_i^l$ for some positive integer $l$.
If $G$ has a cyclic Sylow $p_j$-subgroup of order $p_j^{\alpha_j}$, say $\langle u\rangle$, the $x\sim u\sim y$ is a path, and so $d(x,y)=2$, as desired.
Otherwise, Lemma~\ref{eresult} implies that $p_j=2$ and the Sylow $p_j$-subgroup of $G$ is a generalized quaternion group.
By \eqref{gq} and \eqref{gq-1}, it follows that there exists an element $v$ of order $4$ such that $v\notin \langle y\rangle$.
As a consequence, $x\sim v\sim y$ is a path, and hence $d(x,y)=2$, as desired.
$\qed$

\begin{lemma}\label{lem1}
Let $|G|=p_1^{\alpha_1}p_2^{\alpha_2}\cdots p_t^{\alpha_t}$ where $t\ge 2$. Suppose that for any prime $p_i$, $G$ has a unique subgroup of order $p_i$, where $1\le i \le t$.
If $G$ has an element of order $$p_1^{\alpha_1}p_2^{\alpha_2}\cdots p_{k-1}^{\alpha_{k-1}}
  p_{k}^{\beta_{k}}p_{k+1}^{\alpha_{k+1}}\cdots p_t^{\alpha_t},$$ where $1\le k \le t$ and $1\le \beta_{k}< \alpha_{k}$, then
$\diam(\overline{P(G)^{\ast}})\le 3$.
\end{lemma}
\proof
Let $x,y$ be two distinct vertices such that $x$ and $y$ are non-adjacent in $\overline{P(G)^{\ast}}$.
Without loss of generality, let $\langle x\rangle\subseteq \langle y\rangle$. If $|G|/o(y)\ne p_i^{\beta_i}$ where $1\le \beta_i < \alpha_i$, then it follows from Lemma~\ref{lem-0} that $d(x,y)=2$. Therefore, in the following, we may assume that
$$
o(y)=p_1^{\alpha_1}p_2^{\alpha_2}\cdots p_{k-1}^{\alpha_{k-1}}
  p_{k}^{\beta_{k}}p_{k+1}^{\alpha_{k+1}}\cdots p_t^{\alpha_t},
$$
where $1\le k \le t$ and $1\le \beta_{k}< \alpha_{k}$.
Let $a\in G$ with $o(a)=p_k$, and let
$b\in G$ such that $o(b)$ is a prime with $o(b)\ne p_k$.
Suppose that for some prime $p\ne p_k$, $G$ has at least two  distinct Sylow $p$-subgroups that are cyclic. Let $\langle u\rangle$ be a Sylow $p$-subgroup of $G$ with $\langle u\rangle \nsubseteq \langle y\rangle$.
If $o(x)=p^l$, then $x\sim a\sim u \sim y$ is a path, and so $d(x,y)\le 3$; If not, $x\sim u \sim y$ is a path, and so $d(x,y)=2$.

Thus, we may assume that for any prime $p\ne p_k$, the Sylow $p$-subgroup of $G$ is unique.
Suppose that the Sylow $p_k$-subgroups of $G$ are cyclic.
Let $\langle c\rangle$ be a Sylow $p_k$-subgroup of $G$.
If $o(x)=p_k^l$ where $l$ is a positive integer, then
$x\sim b \sim c \sim y$ is a path, and so $d(x,y)\le 3$.
For the remaining cases,
$x \sim c \sim y$ is a path, which implies $d(x,y)\le 2$,
as desired.

Suppose now that a Sylow $p_k$-subgroup $P_k$ of $G$ is non-cyclic.
Then $P_k$ is isomorphic to a generalized
quaternion group of order at least $8$ by Lemma~\ref{eresult}.
Then we can choose an element $h$ of order $4$ such that $h\notin \langle y\rangle$ by \eqref{gq} and \eqref{gq-1}.
It follows that $x\sim b\sim h \sim y$ is a path for $o(x)=2^l$, and $x\sim h \sim y$ is a path for other cases, as desired.
$\qed$

\begin{lemma}\label{lem2}
Let $G$ be a $\Psi$-group. Then $\diam(\overline{P(G)^{\ast}})=3$.
\end{lemma}
\proof
By the definition of a $\Psi$-group, let
$|G|=p_1^{\alpha_1}p_2^{\alpha_2}\cdots p_t^{\alpha_t}$ with $t\ge 2$, and let $y\in G$ with
$$
o(y)=p_1^{\alpha_1}p_2^{\alpha_2}\cdots p_{k-1}^{\alpha_{k-1}}
  p_{k}^{\beta_{k}}p_{k+1}^{\alpha_{k+1}}\cdots p_t^{\alpha_t},
$$
where $1\le k \le t$ and $1\le \beta_{k}< \alpha_{k}$.
Take $x\in \langle y\rangle$ with $o(x)=p_k$.
By Lemma~\ref{lem1}, we can deduce $d(x,y)\le 3$.
Clearly, we have $d(x,y)\ge 2$.
Now, it suffices to prove $d(x,y)\ne 2$.

Suppose for the sake of contradiction that $d(x,y)=2$. Let $z$ be an element such that $x\sim z \sim y$ is a path.

\medskip
\noindent {\bf Case 1.} $p_k\nmid o(z)$.
\medskip

Since for any prime $p\ne p_k$, the Sylow $p$-subgroup of $G$ is unique, it follows that every Sylow subgroup of $\langle z\rangle$ is contained in $\langle y\rangle$.
Hence, we have $\langle z\rangle \subseteq \langle y\rangle$, and so $z$ and $y$ are non-adjacent in $\overline{P(G)^{\ast}}$, a contradiction.

\medskip
\noindent {\bf Case 2.} $p_k\mid o(z)$.
\medskip

Note that for any prime divisor $p$ of $|G|$, $G$ has a unique subgroup of order $p$. We conclude that $\langle x\rangle\subseteq \langle z\rangle$, and it follows that $z$ and $x$ are non-adjacent in $\overline{P(G)^{\ast}}$, a contradiction.
$\qed$

\medskip
We are now ready to prove Theorem~\ref{mainthm}.
\medskip

\noindent {\em Proof of Theorem~{\rm\ref{mainthm}}.}
By Lemmas~\ref{p-group} and \ref{lem2}, we only need to prove that if $G$ is not $p$-group and is not a $\Psi$-group, then
$\diam(\overline{P(G)^{\ast}})=2$.

Suppose now that $G$ is not a $\Psi$-group with
$|G|=p_1^{\alpha_1}p_2^{\alpha_2}\cdots p_t^{\alpha_t}$, where $t\ge 2$. Clearly, $\diam(\overline{P(G)^{\ast}})\ge 2$.
In view of Lemma~\ref{p-subgroup}, for any prime divisor $p$, we may assume that $G$ has a unique subgroup of order $p$.
Let now $x,y$ be two distinct vertices such that $x$ and $y$ are non-adjacent in $\overline{P(G)^{\ast}}$.
Without loss of generality, we say $\langle x\rangle\subseteq \langle y\rangle$. It follows from Lemma~\ref{lem-0} that we may assume that
$$
o(y)=p_1^{\alpha_1}p_2^{\alpha_2}\cdots p_{k-1}^{\alpha_{k-1}}
  p_{k}^{\beta_{k}}p_{k+1}^{\alpha_{k+1}}\cdots p_t^{\alpha_t},
$$
where $1\le k \le t$ and $1\le \beta_{k}< \alpha_{k}$.
Since $G$ is not a $\Psi$-group, we deduce that there exists at least one prime $p\ne p_k$ such that $G$ has two distinct
Sylow $p$-subgroups.
Consequently, we can
choose a Sylow $p$-subgroup $P$ of $G$ such that $P\nsubseteq \langle y\rangle$.
Note that $|P|$ is not a prime. Also, Lemma~\ref{eresult} implies $P$ is either cyclic or generalized quaternion.
It follows that there exists an element $a$ such that its order is a power of $p$ and $a\notin \langle y\rangle$.

\medskip
\noindent {\bf Case 1.} $p_k\ne 2$.
\medskip

Let $\langle w\rangle$ be a Sylow $p_k$-subgroup of $G$.
If $o(x)=p^l$, then $x\sim w \sim y$ is a path, as desired.
Otherwise, there exists a prime $p'\ne p$ such that $p'\mid o(x)$. This forces that $x\sim a \sim y$ is a path, as desired.

\medskip
\noindent {\bf Case 2.} $p_k=2$.
\medskip

If every Sylow $p_k$-subgroup of $G$ is cyclic, then it is similar to Case 1, the desired result follows.
Assume now that every Sylow $p_k$-subgroup of $G$ is generalized quaternion. Then by \eqref{gq} and \eqref{gq-1},
we can choose an element $u$ of order $4$ such that $u\notin \langle y\rangle$.
If $o(x)=p^l$, then $x\sim u \sim y$ is a path, as desired.
Otherwise, there exists a prime $p'\ne p$ such that $p'\mid o(x)$, and so $x\sim a \sim y$ is a path, as desired.

We conclude $\diam(\overline{P(G)^{\ast}})=2$, the proof is now complete.
$\qed$

\section{Proof of Theorem~\ref{mainthm1}}

In this section, we will prove Theorem~\ref{mainthm1}.

\begin{lemma}\label{e-lem1}
Let $G$ be a finite non-cyclic group. Then
$\diam(\overline{P_e(G)^{\ast}})\le 3$.
\end{lemma}
\proof
If $\overline{P_e(G)^{\ast}}$ is complete, then $\diam(\overline{P_e(G)^{\ast}})=1$, as desired. Thus, in the following we may assume that $\overline{P_e(G)^{\ast}}$ is not complete.
Let $x,y$ be two distinct non-adjacent vertices in $\overline{P_e(G)^{\ast}}$. Then $\langle x,y\rangle$ is cyclic. Suppose that one of $x$ and $y$ can generate a maximal cyclic subgroup of $G$. Without loss of generality, say $\langle y\rangle\in \mathcal{M}(G)$. Since $\langle x,y\rangle$ is cyclic, we have $x\in \langle y\rangle$.
Note that $x\notin \Cyc(G)$. It follows that there exists $\langle z\rangle\in \mathcal{M}(G)$ such that $x\notin \langle z\rangle$. Clearly, $\langle z\rangle \ne \langle y\rangle$. As a result, $x\sim z \sim y$ is a path, and so $d(x,y)=2$.

Thus, in the following, we may assume that $\langle x\rangle\notin \mathcal{M}(G)$ and $\langle y\rangle\notin \mathcal{M}(G)$. Note that $\mathcal{M}_x\ne \mathcal{M}(G)$, $\mathcal{M}_y\ne \mathcal{M}(G)$ and $\mathcal{M}_{\{x,y\}}\ne \mathcal{M}(G)$. If there exists $\langle u\rangle \in \mathcal{M}(G)\setminus \mathcal{M}_{\{x,y\}}$ such that $\{x,y\}\cap \langle u\rangle=\emptyset$, then it is clear that $x \sim u \sim y$ is a path, and so $d(x,y)=2$, as desired.

Thus, now we may assume that for any $M\in \mathcal{M}(G)\setminus \mathcal{M}_{\{x,y\}}$, either $x\in M$ or $y\in M$. Namely, in this case, $G$ is a $\Phi$-group.
Taking now $\langle a\rangle\in \mathcal{M}(G)\setminus \mathcal{M}_{\{x,y\}}$, without loss of generality, we let $x\in \langle a\rangle$. Hence $d(y,a)=1$.
Since $\mathcal{M}_x\ne \mathcal{M}(G)$, there exists $\langle b\rangle\in \mathcal{M}(G)$ such that $x\notin \langle b\rangle$. As a consequence, we have $d(b,x)=1$. Also, it is straightforward that $\langle a\rangle \ne \langle b\rangle$, and so $d(a,b)=1$, which implies that $x \sim b \sim a \sim y$ is a path. Therefore, we have $d(x,y)\le 3$, as desired.
$\qed$

\medskip

By the proof of Lemma~\ref{e-lem1}, the next result is valid.

\begin{lemma}\label{e-lem1-1}
Let $G$ be a non-$\Phi$-group. Then
$\diam(\overline{P_e(G)^{\ast}})\le 2$.
\end{lemma}

\begin{lemma}\label{e-lem2}
Let $G$ be a  $\Phi$-group. Then
$\diam(\overline{P_e(G)^{\ast}})=3$.
\end{lemma}
\proof
Let $x,y\in V(\overline{P_e(G)^{\ast}})$ such that the three conditions in the definition of a $\Phi$-group hold.
Then $x$ and $y$ are non-adjacent in $\overline{P_e(G)^{\ast}}$.
By Lemma~\ref{e-lem1}, it suffices to prove $d(x,y)=3$.

Suppose for a contradiction that $d(x,y)=2$. Let $x\sim z\sim y$ is a path in $\overline{P_e(G)^{\ast}}$.  Then both
$\langle x,z\rangle$ and $\langle y,z\rangle$ are non-cyclic.
Let now $\langle w\rangle\in \mathcal{M}(G)$ with $z\in \langle w\rangle$. It follows that $\{x,y\}\cap \langle w\rangle=\emptyset$, and so $\langle w\rangle\in \mathcal{M}(G)\setminus \mathcal{M}_{\{x,y\}}$,
this contradicts the condition that for any $M\in \mathcal{M}(G)\setminus \mathcal{M}_{\{x,y\}}$, either $x\in M$ or $y\in M$.
$\qed$

\medskip
We are now ready to prove Theorem~\ref{mainthm1}.
\medskip

\noindent {\em Proof of Theorem~{\rm\ref{mainthm1}}.}
For a non-cyclic group $G$, it is clear that $\diam(\overline{P_e(G)^{\ast}})=1$ if and only if $G$ is an elementary abelian $2$-group.
Thus, Theorem~\ref{mainthm1} follows from
Lemmas~\ref{e-lem1}, \ref{e-lem1-1} and \ref{e-lem2}.
$\qed$

\bigskip

\noindent \textbf{Acknowledgements}~~
Ma is supported by NSFC (No. 12326333) and the
Shaanxi Fundamental Science Research Project for Mathematics and Physics (Grant No. 22JSQ024). Wang is supported by the National Key R\&D Program of China (No. 2020YFA0712900) and NSFC (Nos. 12071039 and 12131011).

\end{document}